\documentclass[usenatbib,12pt,letterpaper, notitlepage]{article}
\usepackage[margin=0.9in]{geometry}
\usepackage{graphicx}
\usepackage{amssymb}
\usepackage{amsmath}
\usepackage{amsbsy}
\usepackage{epstopdf}
\usepackage{pdflscape}
\usepackage{setspace}
\usepackage{nicefrac}
\usepackage{multirow}
\usepackage{booktabs}
\usepackage{natbib}
\usepackage{enumerate}
\usepackage{makeidx}  
\usepackage[titletoc]{appendix}
\usepackage{fancyhdr}
\usepackage[parfill]{parskip}	
\usepackage[latin1]{inputenc}
\usepackage{geometry}
\usepackage{graphicx}
\usepackage{epstopdf}
\usepackage{amsmath}
\usepackage{amsfonts}
\usepackage{amssymb}
\usepackage{color}
\usepackage{bm}
\usepackage{float}
\usepackage{setspace}		
\onehalfspace
\usepackage{hyperref}
\usepackage{subfigure}

\allowdisplaybreaks

\begin{document}

\title{
A Broadened Approach to\\ Improved Estimation in Survey Sampling
}
\author{Kyle Vincent and Christopher S. Henry \footnote{Currency Department,
Bank of Canada, 234 Wellington Street,
Ottawa, Ontario, CANADA, K1A 0G9
\textit{email}: kvincent@bankofcanada.ca, chenry@bankofcanada.ca}}

\date{}

\maketitle

\begin{abstract}
\begin{small}
We outline a new framework for design-based inference in a survey sampling setting; we consider how estimation changes when a strategy (that is, a pairing of a sampling design and estimator) is first selected at random amongst a set of candidate strategies. Extending on the traditional approach in survey sampling, we redefine the observed data to now include a mention of the strategy used. The minimal sufficient statistic for the population parameter vector is derived for this setup; data reduction now involves considering hypothetical sampling designs that could have given rise to the sample along with their accompanying estimator(s). Implications on Rao-Blackwellization are investigated. Two studies are detailed to highlight how the new inference procedure can offer alternative views to, and circumvent some challenges encountered with, the traditional approach to inference.
\end{small}
\end{abstract}

\renewcommand{\abstractname}{Acknowledgements}
\begin{abstract}
All views expressed in this manuscript are solely those of the authors and should not be attributed to the Bank of Canada.
\end{abstract}

\clearpage

\section{Introduction}

Traditional design-based inference in survey studies, also known as the fixed population approach, is based on an ordered set of observed sample unit labels together with their corresponding responses of interest. Well-known sources, like those of \cite{Cassel1977} and \cite{Thompson1996}, will denote the observed sample data as $d_0=((i,y_i):i\ \epsilon\ \underline{s})$, where $\underline{s}$ is the set of units in the original order they were selected and $y_i$ is the response of unit $i$. One may notice that there is no mention of the sampling design in the observed data, nor is there any mention of details regarding the estimation procedure. In this manuscript we show that if the strategy is chosen with some probability from amongst a set of candidate strategies then this results in a minimal sufficient statistic based on more than just the population-level data. The new setup to inference results in an approach that has potential to overcome some limitations that are present with the traditional approach, primarily because Rao-Blackwellization now entails averaging over estimators corresponding with multiple strategies.

We consider two cases where the new approach can lend itself well for estimation. First, \cite{Cassel1977} reviewed the admissibility properties of strategies and noted how no strategy is guaranteed to always perform better than any other strategy. Consequently, at times it may be unclear as to which strategy should serve best for a study. Further, some well-known strategies can give rise to highly-skewed sampling distributions for commonly used estimators, like when one is sampling from highly uneven populations with probability proportional to size depending on auxiliary information and the Hansen-Hurwitz estimator is used. Upon selecting such samples the analyst may feel tempted or inclined to make inference based on a different strategy. With our setup, we show how the analyst can exploit the corresponding Rao-Blackwell features to make this possible, which in turn can dampen the impact of designs on estimators sensitive to such possibly unrepresentative samples. Hence, with our method one can properly base inference on a conservative approach through taking a composite average of several strategies to avoid committing to one potentially weak strategy.

Second, we consider how the new setup can aid in a retrospective style of analysis. We consider situations when only a limited amount of information is available on the actual sampling design that was used, possibly because the procedure is beyond the full control of the analyst. We show how the analyst can utilize such information in the analysis with the aid of the new setup we present. The method is interesting in that it parallels the Bayesian paradigm as it allows the analyst to draw on prior knowledge of the design used to select the sample, as reflected upon in the discussion section of this manuscript.

Our goal is to offer a new and alternative view to design-based inference with ideas as to where the approach may be practical in an empirical setting. We consider empirical data sets and simulation studies to justify our claims and encourage future work based on the new setup. The article is organized as follows. Section 2 introduces the setup with mathematical notation. Section 3 provides the minimal sufficiency result and a proof that the minimal sufficient statistic is not complete. Section 4 gives details and formulas oriented about the improved estimation procedure, which includes the variance expression and a corresponding expression for an estimate of the variance of the improved estimator. Section 5 provides several ideas with simulation results to demonstrate how the new setup is practical for some specific cases. Finally, Section 6 discusses the implications of the novel inference procedure outlined in this manuscript.

\section{Framework and Mathematical Notation}

Define $U=\{1,2,...,N$\} to be the population unit labels where $N$ is the size of the population. Define $\underline{\theta}=(y_1,y_2,...,y_N)$ to be the population parameter vector; $y_i$ is the response of unit $i$, $\underline{\theta}\ \epsilon\ \Theta$ where $\Theta=\mathcal{R}^{N}$, and $\mathcal{R}$ is the set of all real numbers.

Define $\mathcal{D}$ to be the set of candidate sampling designs under consideration for the data collection aspect of the study. For ease of presentation we will assume that there are a finite number of sampling designs. Hence, the set of candidate designs can be ordered such that $\mathcal{D}=(\delta_1,\delta_2,...,\delta_T)$ where $\delta_k$ refers to sampling design $k$, $k=1,2,...,T$. For example, $\delta_1$ may refer to the simple random sampling design and $\delta_2$ may refer to a probability proportional to size sampling design.

Suppose $\phi=\phi(\underline{\theta})$ is a population quantity to be estimated. Define $\Phi_k$ to be the set of candidate estimators that accompany sampling design $k$. For ease of presentation we will assume that there are a finite number of estimators accompanying $\Phi_k$. Hence, the corresponding set of candidate estimators can be ordered such that $\Phi_k=(\hat{\phi}_{k,1},\hat{\phi}_{k,2},...,\hat{\phi}_{k,V_k})$ where $\hat{\phi}_{k,j}$ refers to estimator $j$ that accompanies sampling design $k$, $j=1,2,...,V_k$. For example, if $\delta_1$ is the simple random sampling design then $\hat{\phi}_{1,1}$ may refer to the sample mean estimator and $\hat{\phi}_{1,2}$ may refer to the sample median estimator.

Define $\Lambda_{k,j}$ to be the strategy that corresponds with sampling design $\delta_k$ and estimator $\hat{\phi}_{k,j}$. Let $\underline{P}_{\Lambda}=(P(\Lambda_{k,j}):k=1,2,...,T, j=1,2,...,V_k)$ be the probability vector with $P(\Lambda_{k,j})$ referring to the probability that strategy $\Lambda_{k,j}$ is selected for data collection and inference.

Define $\underline{s}$ to be the vector of sample unit labels displayed in the order they were selected (with possible repeat selections). Define $P(\underline{s}|\underline{y}_{\underline{s}},\delta_k)$ to be the probability that $\underline{s}$ is selected under sampling design $k$ where $\underline{y}_{\underline{s}}$ is the corresponding vector of observed responses for all individuals selected for the original sample $\underline{s}$.

The original data that is observed is $d_0=(\Lambda_{k,j},(i,y_i):i\ \epsilon\ \underline{s})$. Define $r_d$ to be the reduction function that maps the original data to the reduced data $d_R=\{(i,y_i):i\ \epsilon\ s\}$ (that is $r_d(d_0)=d_R$) where $s$ is the unordered set of unit labels. Define the reduced data ordered from smallest to largest label in the sample to be $d=(s,\underline{y}_s)$. As is common in the framework of design-based inference (see \cite{Thompson1996} for a comprehensive review of the design-based setup to inference in survey sampling), we define a parameter vector $\underline{\theta}$ to be \textit{consistent} with $d$ if the $i^{\text{th}}$ component of $\underline{y}_s$ is equal to the $s_{i}^{\text{th}}$ component of $\underline{\theta}$, say. We define $\Theta_{d}$ to be the set of all $\underline{\theta}$ that are consistent with $d$. Notice that, since consistency depends only on the values of the distinct units in the sample, $\Theta_{d}=\Theta_{d_R}$.

\section{Minimal Sufficiency Result and Incompleteness}

In the traditional setup the minimal sufficient statistic for $\underline{\theta}=(y_1,y_2,...,y_N)$ has been shown to be $d_R=\{(i,y_i):i\ \epsilon\ s\}$ (\cite{Pathak1964, Godambe1966, Basu1969}). \cite{Cassel1977} later showed that the minimal sufficient is not complete. In this section we show that the minimal sufficient statistic in our inferential setup is still the aforementioned reduced data and that it is also not complete.

\bigskip

\noindent \textbf{Theorem:} Suppose that $\underline{P}_{\Lambda}$ does not depend on $\underline{\theta}$ and that all candidate sampling designs are of the adaptive and/or conventional type (see \cite{Thompson1996} for a discussion on such designs). That is, for all $\underline{s}$ and $\Lambda_{k,j}$, $P(\underline{s}|\underline{y},\Lambda_{k,j})=P(\underline{s}|\underline{y}_{\underline{s}},\delta_k)(=P(\underline{s}|\delta_k)$ if the design is conventional). The random variable $D_R$ is the minimal sufficient statistic for $\underline{\theta}$.

\bigskip

\noindent \textbf{Proof}: Choose any possible $d_0$ and consider $\Theta_{d_R}$. Suppose that strategy $\Lambda_{k,j}$ is the strategy selected for the data collection and inference aspects of the study, for some arbitrary $k$ and $j$. Now,
\begin{align}\label{Sufficiency Theorem}
P_{\underline{\theta}}(D_0=d_0)&=P(\Lambda_{k,j})P(\underline{s}|\underline{y}_{\underline{s}},\delta_k)I[\underline{\theta}\ \epsilon\ \Theta_{d_R}]\\\notag
&=L_{\underline{\theta}}(D_0=d_0).
\end{align}

As the likelihood can be separated into two components, one which does not depend on $\underline{\theta}$ and the other which depends on $\underline{\theta}$ only through the reduced data, by the Neyman Factorization Theorem we can conclude that $D_R$ is a sufficient statistic for $\underline{\theta}$.

To show the minimality of the claim, we will make use of the theorem which gives the following. If $\underline{X}$ is a sample with probability mass/density function $f(\underline{X}|\underline{\theta})$ and $T(\underline{X})$ is a function of $\underline{X}$ such that, for any two sample points $\underline{x}$ and $\underline{x}^*$, the function $k(f(\underline{x}),f(\underline{x}^*))>0$ is constant as a function of $\underline{\theta}$ if and only if $T(\underline{x})=T(\underline{x}^*)$, then $T(\underline{X})$ is the minimal sufficient statistic for $\underline{\theta}$ \citep{Arnold1990, Casella2002}.

Take any $d_0=(\Lambda_{k,j},(i,y_i):i\ \epsilon\ \underline{s})$ and $d_0^*=(\Lambda_{k^*,j^*},(i,y_i):i\ \epsilon\ \underline{s}^*)$ where $P(d_0),P(d_0^*)>0$. Suppose that $P_{\underline{\theta}}(d_0)=k(d_0,d_0^*)P_{\underline{\theta}}(d_0^*)$ where $k$ is independent of $\underline{\theta}$. We can re-express this as
\begin{align}
P(\Lambda_{k,j})P(\underline{s}|\underline{y}_{\underline{s}},\delta_k)I[\underline{\theta}\epsilon\Theta_{d_R}]=k(d_0,d_0^*)P(\Lambda_{k^*,j^*})P(\underline{s}^*|\underline{y}_{\underline{s}^*},\delta_{k^*})I[\underline{\theta}\epsilon\Theta_{d_R^*}]
\end{align}

\noindent where $d_R=\{(i,y_i):i\ \epsilon\ s\}, d_R^*=\{(i,y_i):i\ \epsilon\ s^*\}$. As
\begin{align}
P(\Lambda_{k,j}),P(\Lambda_{k^*,j^*}),P(\underline{s}|\underline{y}_{\underline{s}},\delta_k),P(\underline{s}^*|\underline{y}_{\underline{s}^*},\delta_k^*)>0,
\end{align}
\noindent the indicators must take on a value of zero or one at the same time. Hence, $I[\underline{\theta}\epsilon\Theta_{d_R}]=I[\underline{\theta}\epsilon\Theta_{d_R^*}]$ for all $\underline{\theta}$ and so $\Theta_{d_R}=\Theta_{d_R^*}$, which implies that $d_R=d_R^*$. Therefore, $D_R$ is the minimal sufficient statistic.

$\Box$

\bigskip

\noindent \textbf{Theorem}: The random variable $D_R$ is not complete.

\bigskip

\noindent \textbf{Proof}: By way of counterexample, consider the case where only one strategy is considered. Then this is the usual case in the classical fixed population approach since mentioning the strategy in the observed data is redundant. Hence, by \cite{Cassel1977} it must be that $D_R$ is not complete.

$\Box$

\section{Estimation}

In this section we derive the expression for the improved point estimators and an estimate of the variance of this estimator. We conclude the section by highlighting how the suggested method of estimation for the variance of the improved estimator avoids the need for estimation of covariances.


\subsection{Point estimation}

Suppose $\phi=\phi(\underline{\theta})$ is a population quantity to be estimated. Define $\mathcal{S}$ to be the set of all subsets of $U$ that have a positive probability of being obtained with at least one $\delta_k\ \epsilon\ \mathcal{D}$. For ease of presentation, we will define
\begin{align}
\hat{\phi}_{k,j}(s)=\frac{\sum\limits_{\underline{s}\epsilon\mathcal{S}_s}\hat{\phi}_{k,j}(\underline{s})P(\underline{s}|\underline{y}_{\underline{s}},\delta_k)}
{\sum\limits_{\underline{s}\epsilon\mathcal{S}_s}P(\underline{s}|\underline{y}_{\underline{s}},\delta_k)}
\end{align}
where $\mathcal{S}_s$ is the set of all full samples $\underline{s}$ such that the corresponding reduced set is $s$. Notice that this is the usual Rao-Blackwellization scheme that is used in traditional survey sampling inference.

Recall that we defined $\Lambda_{k,j}$ to be the strategy that depends on sampling design $\delta_k$ and estimator $\hat{\phi}_{k,j}$. With respect to any specific sample $s$, in the event that $s$ cannot be obtained with $\delta_k$, that is $P(s|\underline{y}_s,\delta_k)=0$, we will define $\hat{\phi}_{k,j}(s)$ to be zero.

Let $\hat{\phi}(s)$ denote the preliminary estimate of $\phi$ (preliminary in the sense that improvement has not been made with respect to averaging over estimators corresponding with multiple strategies). The Rao-Blackwellized version of $\hat{\phi}(s)$ is
\begin{align}\label{RB expression}
\hat{\phi}_{RB}(s)&=E[\hat{\phi}(s)|d_R]\\\notag
&=\frac{\sum\limits_{k}\sum\limits_{j}(\hat{\phi}_{k,j}(s)P(\Lambda_{k,j})P(s|\underline{y}_s,\delta_k))}{\sum\limits_{k}\sum\limits_j(P(\Lambda_{k,j})P(s|\underline{y}_s,\delta_k))}\\\notag
&=\frac{\sum\limits_{k}\sum\limits_{j}(\hat{\phi}_{k,j}(s)P(\Lambda_{k,j})P(s|\underline{y}_s,\delta_k))}{P(s|\underline{y}_s)}.
\end{align}

\subsection{Variance estimation}

The variance of the preliminary estimate $\hat{\phi}$ is
\begin{align}
\text{var}(\hat{\phi})&=\sum\limits_{s\epsilon \mathcal{S}}\sum\limits_{k}\sum\limits_{j}((\hat{\phi}_{k,j}(s)-E[\hat{\phi}])^2P(\Lambda_{k,j})P(s|\underline{y}_s,\delta_k)).
\end{align}

The variance of the improved estimator can be determined with the decomposition of variances as follows,
\begin{align}
&\text{var}(\hat{\phi}_{RB})=\text{var}(\hat{\phi})-E[\text{var}(\hat{\phi}|d_R)]\\\notag
=&\sum\limits_{s\epsilon\mathcal{S}}\sum\limits_{k}\sum\limits_{j}((\hat{\phi}_{k,j}(s)-E[\hat{\phi}])^2P(\Lambda_{k,j})P(s|\underline{y}_s,\delta_k))\\\notag
&-\sum\limits_{s\epsilon \mathcal{S}}\sum\limits_{k}\sum\limits_{j}((\hat{\phi}_{k,j}(s)-\hat{\phi}_{RB}(s))^2P(\Lambda_{k,j})P(s|\underline{y}_s,\delta_k))\\\notag
=&\sum\limits_{s\epsilon\mathcal{S}}\sum\limits_{k}\sum\limits_{j}((\hat{\phi}_{RB}(s)-E[\hat{\phi}])^2P(\Lambda_{k,j})P(s|\underline{y}_s,\delta_k))\\\notag
=&\sum\limits_{s\epsilon\mathcal{S}}((\hat{\phi}_{RB}(s)-E[\hat{\phi}])^2P(s|\underline{y}_s)).
\end{align}

To estimate the variance of the improved estimate, the following estimator can be used. First, an estimate of $\text{var}(\hat{\phi})$ is
\begin{align}\label{Var estimate first}
E[\hat{\text{var}}(\hat{\phi}(s))|d_R]=\sum\limits_{k}\sum\limits_{j}  (\hat{\text{var}}(\hat\phi_{k,j}(s))P(\Lambda_{k,j})P(s|\underline{y}_s,\delta_k))/P(s|\underline{y}_s)
\end{align}

where $\hat{\text{var}}(\hat{\phi}_{k,j}(s))$ is the corresponding estimate of the variance of the estimate of $\hat{\phi}_{k,j}(s)$. Notice that this is the Rao-Blackwellized estimate of the preliminary estimate of the variance of $\hat{\phi}_{k,j}(s)$.

Second, an estimate of $E[\text{var}(\hat{\phi}|d_R)]$ is
\begin{align}\label{Var estimate second}
\text{var}(\hat{\phi}(s)|d_R)=\sum\limits_k\sum\limits_{j}((\hat{\phi}_{k,j}(s)-\hat{\phi}_{RB}(s))^2P(\Lambda_{k,j})P(s|\underline{y}_s,\delta_k))/P(s|\underline{y}_s).
\end{align}

Finally, the estimate of $\hat{\text{var}}(\hat{\phi}_{RB})=E[\hat{\text{var}}(\hat{\phi}(s))|d_R]-\text{var}(\hat{\phi}(s)|d_R)$ can be used as an estimate of $\text{var}(\hat{\phi}_{RB})$. Furthermore, if $\hat{\text{var}}(\hat{\phi}_{k,j}(s))$ is an unbiased estimate of $\text{var}(\hat{\phi}_{k,j}(s))$ for all $k$ and $j$, then this estimator is unbiased. We make two remarks here. First, a direct attempt to determining the variance of the improved estimator could be challenging as this would require determining the covariance of estimators corresponding with multiple strategies. In contrast, the aforementioned estimator has the advantage in that it avoids covariance measures of estimators. Second, in some cases this estimator will give negative estimates. In such a case, a conservative approach is to take $\hat{\text{var}}(\hat{\phi}_{RB})=E[\hat{\text{var}}(\hat{\phi}(s))|d_R]$. However, this is still an outstanding issue and will require future attention.



\section{Simulation Studies}

In this section we demonstrate how the new inference procedure can assist in various aspects of a survey sampling study. Results from two simulation studies are provided.

\subsection{Simulation study 1}\label{simstudy1}

The primary objective of this simulation study is to demonstrate how the new inference procedure can allow one to average over estimates corresponding with several strategies to avoid committing to, and possibly obtain a more efficient estimator than that based on, any one of the candidate strategies. We consider the `influenza' data set in R \citep{samplingbook} where a count of inhabitants and influenza cases within 424 districts of Germany in 2007 are provided. We aim to estimate the total number of cases of influenza. Two strategies are considered. The first strategy, denoted $\Lambda_{1,1}$, entails the use of a simple random sampling with replacement (SRSWR) design and the full sample mean. The second strategy, denoted $\Lambda_{2,1}$, entails the use of a probability proportional to size with replacement (PPSWR) sampling design and a Hansen-Hurwitz (HH) estimator where the count of inhabitants serves as the auxiliary information.  As we are interested in how the new setup can improve on estimation, Rao-Blackwellization is restricted to directly averaging over the aforementioned estimators corresponding with the two strategies and not with respect to the unique sample elements\footnote{Though sampling is carried out with replacement, we consider small sample sizes relative to the population size and hence one may argue that the use of the Hansen-Hurwitz estimator is practical for such a case. Further, one can avoid the tedious computations required for inclusion probabilities typically required for without-replacement estimators.}.

\begin{figure}[H]
\begin{center}
\includegraphics [width=6.5in]{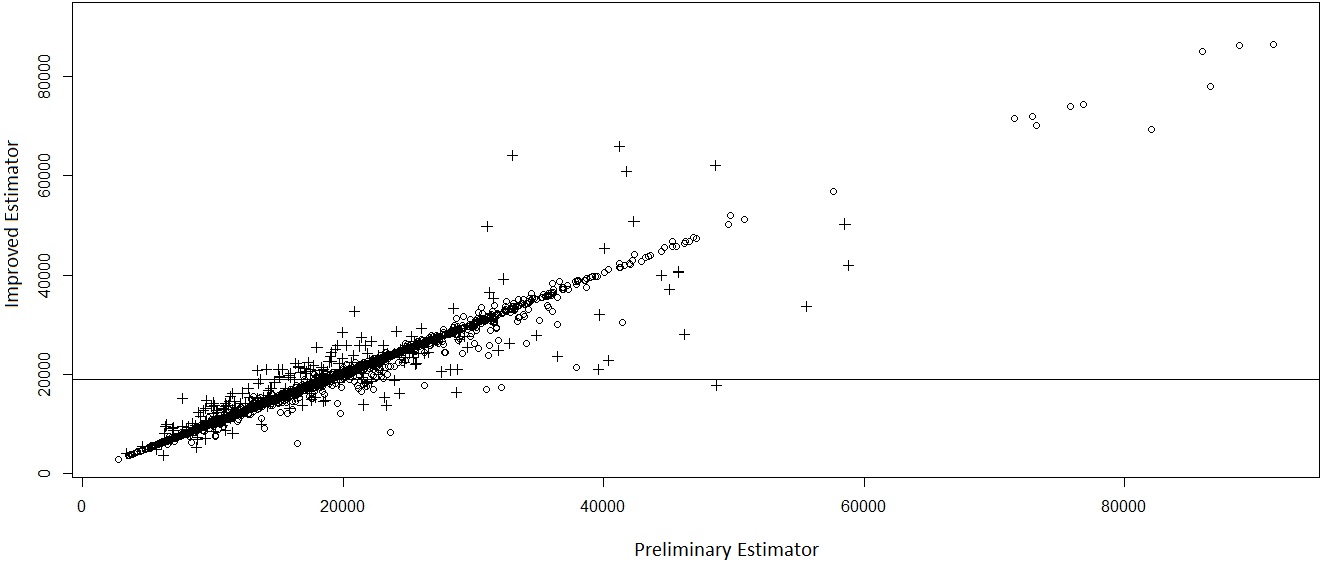}
\end{center}
\caption{A plot of 2000 randomly selected points from the simulation study.  Preliminary estimates correspond with the x-axis and improved estimates correspond with the y-axis. The crosses correspond with the use of $\Lambda_{1,1}$ and circles correspond with the use of $\Lambda_{2,1}$. The vertical and horizontal lines correspond with the true population total.}
\label{results}
\end{figure}

$P(\Lambda_{1,1})$ is set to $0.10$ and $P(\Lambda_{2,1})$ is set to $0.90$ to reflect a level of certainty in the use of a PPSWR design. Samples of size five are obtained and a large number of samples are drawn in order to eliminate any Monte carlo error. Table \ref{Simulation study 1} provides the standardized variance scores and corresponding coverage rates based on nominal 95\% confidence intervals via the Central Limit Theorem (CLT).
\begin{table}[H]
\centering
\caption{Simulation study based on the influenza data set in R. The sample size is five and $P(\Lambda_{1,1})=0.10, P(\Lambda_{2,1})=0.90$. Population total is 18,900. All estimators are unbiased. The variance scores are standardized about the Hansen-Hurwitz estimator used with the PPSWR design (that is, those obtained with strategy $\Lambda_{2,1}$).}
\begin{tabular}{l*{6}{r}r}
Strategy               &Standardized variance scores        &Coverage rates\\\hline
$\Lambda_{1,1}$        &1.330                              &0.971\\
$\Lambda_{2,1}$        &1                                  &0.886\\
Preliminary            &1.033                              &0.895\\
Improved               &0.961                              &0.913\\
\hline
\label{Simulation study 1}
\end{tabular}
\end{table}
In this case the improved estimator outperforms each of the individual strategy-based and preliminary estimators. Further, the coverage rates of the improved estimator rest between those based on the two strategies. To further investigate the details of the study, Figure \ref{results} provides a plot of the preliminary by improved estimators for 2000 randomly selected points. The sampling design used to select the sample is reflected upon by the character of the plot. Notice how the improved estimators are typically pulled towards their expectation, and in general there is a dampening effect on the preliminary estimates with Rao-Blackwellization. This is especially evident when such estimators are extreme. In some cases averaging over the two strategies is intuitive since the sample may \textit{appear} as if it were selected under a design corresponding with a strategy not chosen at the initial stage of the study (in particular when extreme preliminary estimates are reported). Note that if considering designs like SRSWR and PPSWR when sample sizes are small, this could very easily be the case.

The simulation study is extended over a range of sample sizes. It is found that with sizes less than or equal to fifteen, the performance of the improved estimator is noticeably better than each of the strategy-based and preliminary estimators. Hence, we can conclude that in some small sample studies the use of the inference procedure we have outlined may serve as a competitive alternative to, and may result in a more efficient than some of, commonly used estimators based on the traditional setup. For sample sizes greater than fifteen the improved estimator is weighted almost entirely towards the strategy used at the preliminary stage. Hence the benefit from using the new inference procedure in these cases is negligible. This is primarily due to the probability of selecting the sample under the design corresponding with the strategy chosen at the initial stage of the study relative to the probability of selecting it with the design corresponding with the strategy that is not chosen.

\subsection{Simulation study 2}\label{simstudy3}

The primary objective of this simulation study is to demonstrate how the new inference procedure can assist in a retrospective analysis. The details of the study are provided below.

We consider a case where the sample is selected via SRSWR but the sample size is unknown; only the unique elements of the sample are recorded.  Features of the effective sample mean are of interest, in particular a variance estimator of the sample mean. We apply our inference setup to obtain an improved estimate of the variance of the effective sample mean estimator.

The effective sample mean is defined to be $\bar{y}_{\nu}=\frac{1}{\nu}\sum\limits_{i\epsilon s}y_i$ where $\nu$ is the size of the unique sample $s$. The variance of the effective sample mean estimator \citep{Pathak1962} is
\begin{align}\label{var true effective}
\text{var}(\bar{y}_{\nu})=\frac{\sum\limits_{j=1}^{N-1}j^{n-1}}{N^n}\frac{1}{N-1}\sum\limits_{i=1}^N(Y_i-\mu)^2
\end{align}

\noindent where $\mu$ is the population mean and $n$ is the number of draws made to select the sample. An unbiased estimator for $\text{var}(\bar{y}_{\nu})$ is
\begin{align}\label{var effective}
\hat{\text{var}}(\bar{y}_{\nu})=\frac{\sum\limits_{j=1}^{N-1}j^{n-1}}{N^n}\frac{1}{\nu-1}\sum\limits_{i\epsilon s}(Y_i-\bar{y}_{\nu})^2.
\end{align}

We shall define $\delta_k$ to be the SRSWR design that depends on a specific final sample size where repeat selections may occur, $n_k$ say. We shall define $\hat{\phi}_{k,1}$ to be the estimator in \eqref{var effective} that depends on the corresponding sample size associated with $\delta_k$. Hence, the improved version of the estimator presented in \eqref{var effective} is of the form
\begin{align}\label{RB expression var}
\hat{\text{var}}(\bar{y}_{\nu})_{RB}&=E[\hat{\text{var}}(\bar{y}_{\nu})_{n_k}(s)|d_R]\\\notag
&=\frac{\sum\limits_{k}(\hat{\text{var}}(\bar{y}_{\nu})_{n_k}(s)P(\Lambda_{n_k})P(s|\delta_k))}{P(s)}.
\end{align}

where $P(\Lambda_{n_k})$ is the probability that the sample size is $n_k$ and $\hat{\text{var}}(\bar{y}_{\nu})_{n_k}(s)$ is the estimate of the variance of $\bar{y}_{\nu}$ when the sample size is $n_k$. Notice that this estimator will depend on evaluating the probability of selecting the unique sample under a SRSWR design for a series of specific final sample sizes. We have not come across an expression in the literature for evaluating such a probability and have therefore taken the liberty of providing its derivation in the Appendix.


In our study we consider the `trees' dataset \citep{datasets}, a data set that provides measurements of the girth, height, and volume of timber in 31 felled black cherry trees. We define $\Lambda_{k,1}$ to be based on a SRSWR design with a final sample size of $5+k-1$ where $k=1,2,...,10$. In each case, $\hat{\phi}_{k,1}$ is the estimate of the variance of the effective sample mean based on the corresponding sample size of $5+k-1$. In the simulation study strategies are chosen such that $P(\Lambda_{k,1})\dot{\propto}\frac{1}{k}$ so that smaller sample sizes are expected.  Table \ref{Simulation study 2} summarizes the results based on a large number of simulation runs. The naive estimator is based on the actual effective sample size (that is, $n$ is treated as $\nu$ so that inference is based on a hypothetical sample that is equivalent to the unique sample). Coverage results are based on nominal 95\% CLT confidence intervals.

\begin{table}[H]
\centering
\caption{Simulation study 2 results. Results based on the girth, height, and volume responses, respectively.}
\begin{tabular}{l*{6}{r}r}
Response        &Estimator           &Relative Bias      &Variance               &Coverage\\\hline
Girth\\\hline
                &Naive               &0.138              &0.783                  &0.904\\
                &Preliminary         &0                  &0.593                  &0.889\\
                &Improved            &0                  &0.563                  &0.890\\
\hline
\hline
Volume\\\hline
                &Naive               &0.131              &13.00                  &0.911\\
                &Preliminary         &0                  &9.94                   &0.894\\
                &Improved            &0                  &9.46                   &0.895\\
\hline
\hline
Height\\\hline
                &Naive               &0.131              &826                    &0.882\\
                &Preliminary         &0                  &626                    &0.866\\
                &Improved            &0                  &601                    &0.868\\
\hline
\label{Simulation study 2}
\end{tabular}
\end{table}

In each case, the new inference procedure provides improved estimates, and therefore tighter confidence bands corresponding with the estimates for the population mean, as well as coverage rates of the population mean that are on par with those found using the preliminary estimator. Hence, we can conclude that the new inferential setup has the potential to lend itself well for obtaining efficient estimators for sampling distributions when only limited information pertaining to the sampling design is available.

\section{Implications}

In this manuscript we have introduced a new approach to making design-based inference in survey studies; the framework we outline rests on incorporating an additional element of randomization via allowing for the strategy to first be selected at random from amongst a set of candidate strategies. For such a case, we have derived the minimal sufficient statistic for the population parameter vector and have shown that it is not complete. We demonstrate how the new approach can be useful for some specific survey sampling studies.



In general, the new inference framework allows the analyst to make a compromise over a set of candidate strategies through the Rao-Blackwellized estimator and is therefore not required to fully commit to one potentially weak strategy. In the first simulation study we demonstrate that studies based on small sample sizes can benefit from basing a final estimator on an average of estimators corresponding with strategies that are comprised of various sampling designs. As we highlighted, this approach works especially well when extreme values are reported for preliminary estimators. Research on how small sample studies that usually encounter such occurrences can benefit from this approach will make for useful future work.

In an empirical setting the analyst may not have full control over the sample selection procedure. Therefore, they may only have a limited amount of information pertaining to the sampling design that was used for data collection. An example motivates the second simulation study; when sampling is carried out via a with replacement sampling design it could be the case that only the unique elements of the sample are provided to the analyst, perhaps due to a lack of communication. In such occurrences one can posit a series of sample sizes and determine Rao-Blackwell-type estimators for population quantities. We further explored this approach for the simulation parameters outlined in the second simulation study. We consistently sampled at a fixed final sample size and posited sample sizes that ranged up to twice the final sample size. We also posited equal probabilities for the sample sizes. We found that coverage rates corresponding with the improved estimator were only negligibly affected while the variance was significantly reduced.

The approach presented in this paper suggests a new method for applying a Bayes like analysis to a sampling setup; consider Expression \eqref{RB expression} and how it is reflective of a Bayes estimator. In the case of when limited knowledge is available on the sampling design used, one can posit a series of (prior) sampling designs with a corresponding (prior) distribution that each design was used and then base inference on the Rao-Blackwell expression. Future work on combinations of strategies that can be used to draw meaningful conclusions in certain situations when the analyst has limited control/knowledge of the sample selection procedure is deserving of future attention.

\clearpage
\bibliographystyle{biom}
\bibliography{MasterReferences}


\clearpage

\appendix
\section{Appendix}\label{appA}

In this appendix we derive the probability of obtaining a specific sample of unique elements where the size of the unique sample is $\nu$, the size of the full sample is $n$, and the population size is $N$.

Suppose that units are selected from a population where $U=\{1,2,...,N\}$ are the population unit labels and where selections are made via the simple random sampling with replacement (SRSWR) design of size $n$. Note that this is equivalent to drawing one unit at a time (where each unit has probability $1/N$ of being selected on each draw), recording the label identification, and placing the unit back in the population. Suppose that $\underline{s}$ is the original ordered sample. Let $s$ be the effective (unique) sample. For example, if $\underline{s}=(2,3,1,1,4,2)$ then $s=\{1,2,3,4\}$. We would like to determine the probability of obtaining $s$, $P(s)$.

Partition $U$ into $s$ and $U\setminus s$. Now, define $X$ to be the event that all $n$ draws are made from $s$, and define $Y$ to be the event that no units in $s$ are missing from the final sample. We would like to determine $P(s)=P(X\cap Y)=P(X)\times P(Y|X).$ Suppose $|s|=\nu$. Clearly,
$P(X)=(\frac{v}{N})^n$. Now, define $A_i$ to be the event that unit $i$ is selected for the final sample given that all draws are made from $s$, $i\ \epsilon\ s$. We will now determine
\begin{align}
&P(Y|X)=1-P(\cup_{i\epsilon s} \overline{A}_i).
\end{align}
By the principle of inclusion-exclusion,
\begin{align}\label{inclusionexclusion}
P(\cup_{i\epsilon s} \overline{A}_i)=\sum\limits_{k=1}^\nu\bigg((-1)^{k-1}\sum\limits_{I\subset s:|I|=k}P(\overline{A}_I)\bigg).
\end{align}
If $|I|=k$, then $P(\overline{A}_I)=\bigg(\frac{v-k}{v}\bigg)^n$  so that expression \eqref{inclusionexclusion} can be rewritten as
\begin{align}
\sum\limits_{k=1}^\nu\bigg((-1)^{k-1}{\nu \choose k}\bigg(\frac{\nu-k}{\nu}\bigg)^n\bigg)
\end{align}
since there are ${\nu \choose k}$ subsets of $s$ that are of size $k$.

Therefore,
\begin{align}
P(s)=\bigg(\frac{v}{N}\bigg)^n\bigg[1-\sum\limits_{k=1}^\nu\bigg((-1)^{k-1}{\nu \choose k}\bigg(\frac{\nu-k}{\nu}\bigg)^n\bigg)\bigg].
\end{align}

\end{document}